\documentclass[11pt]{article}
\usepackage{amsmath,amssymb,amsfonts,amscd}
\usepackage{mathrsfs}
\usepackage{dsfont}
\usepackage{graphicx}
\usepackage{epsfig,psfrag,eepic}
\usepackage{pifont}
\usepackage{enumitem}
\usepackage{latexsym}
\usepackage{color,colordvi}
\usepackage{cite}

\newtheorem{theorem}{Theorem}[section]

\newtheorem{lemma}{Lemma}[section]

\newtheorem{claim}{Claim}[section]

\def\qed{\hfill \ensuremath{\square}}

\def\pf{\noindent {\it Proof. }}

\usepackage[colorlinks=true,linkcolor=blue,citecolor=blue,urlcolor=blue]{hyperref}

\begin{document}
	
	\title{On Perfect Divisibility of Bull-Free Graphs Without Long Paths}
	\author{
		Di Wu$^{1}$\footnote{
		Email: diwu@njit.edu.cn}, \; Junran Yu$^{2}$\footnote{Email: nnuyu@foxmail.com}, \; Xiaowen Zhang$^{3}$\footnote{ Email: xiaowzhang0128@126.com}\\
		\small $^1$Department of Mathematics and Physics\\
		\small Nanjing Institute of Technology, Nanjing 211167, Jiangsu, China\\
		\small $^2$Institute of Mathematics, School of Mathematical Sciences\\
		\small Nanjing Normal University, 1 Wenyuan Road,  Nanjing, 210023,  China\\
		\small $^3$Department of Mathematics\\
		\small East China Normal University, Shanghai, 200241, China
	}
	\date{}
	\maketitle
	
	\begin{abstract}
		A graph $G$ is {\em perfectly divisible} if, for every induced subgraph $H$ of $G$, $V(H)$ can be partitioned into $A$ and $B$ such that $H[A]$ is perfect and $\omega(H[B])<\omega(H)$. Chudnovsky and Sivaraman [J. Graph Theory \textbf{90} (2019) 54-60]  proved that every ($P_5$, bull)-free graph is perfectly divisible, while Chen and Xu [Discrete Appl. Math. \textbf{372} (2025) 298-307] proved the same for ($P_7,C_5$, bull)-free graphs. We extend these results by proving that every ($P_8,C_5$, bull)-free graph is perfectly divisible and that, letting $F$ denote the Gr\"{o}tzsch graph, a ($P_6$, bull)-free graph is perfectly divisible if and only if it is $F$-free.
	\end{abstract}
	
	\begin{flushleft}
	\end{flushleft}
	
	\section{Introduction}\label{introduction}
	
	Throughout this paper, all graphs are finite and simple. We follow \cite{BM08} for undefined notation and terminology. As usual, $P_k$ and $C_k$ denote the path and the cycle on $k$ vertices, respectively. We say that a graph $G$ {\em contains} a graph $H$ if $H$ is isomorphic to an induced subgraph of $G$. A graph $G$ is {\em $H$-free} if it does not contain $H$. More generally, for a family $\mathcal{H}$ of graphs, $G$ is $\mathcal{H}$-free if it contains no member of $\mathcal{H}$ as an induced subgraph.
	
	Let $H_1$ and $H_2$ be two vertex-disjoint graphs. The {\em join} $H_1+H_2$ is the graph with $V(H_1+H_2)=V(H_1)\cup V(H_2)$ and $E(H_1+H_2)=E(H_1)\cup E(H_2)\cup\{xy\mid x\in V(H_1),\,y\in V(H_2)\}.$
	
	A {\em proper coloring} of a graph $G$ is an assignment of colors to the vertices of $G$ such that adjacent vertices receive different colors. The {\em chromatic number} of $G$, denoted by $\chi(G)$, is the minimum number of colors needed for a proper coloring of $G$. A {\em clique} (resp. {\em independent set}) of $G$ is a set of pairwise adjacent (resp. pairwise nonadjacent) vertices. The {\em clique number} of $G$, denoted by $\omega(G)$, is the maximum size of a clique in $G$.
	
	An induced cycle of length $k\ge 4$ is called a {\em hole}, and $k$ is the {\em length} of the hole. A {\em $k$-hole} is a hole of length $k$. A hole is {\em odd} if its length is odd. An {\em antihole} is the complement of a hole. A graph $G$ is said to be {\em perfect} if $\chi(H)=\omega(H)$ for every induced subgraph $H$ of $G$. The Strong Perfect Graph Theorem, proved by Chudnovsky {\em et al.} \cite{CRST06}, states that:
	
	\begin{theorem}\label{perfect}{\em \cite{CRST06}}
		A graph $G$ is perfect if and only if $G$ is (odd hole, odd antihole)-free.
	\end{theorem}
	
\begin{figure}[htbp]
		\begin{center}
			\includegraphics[width=9cm]{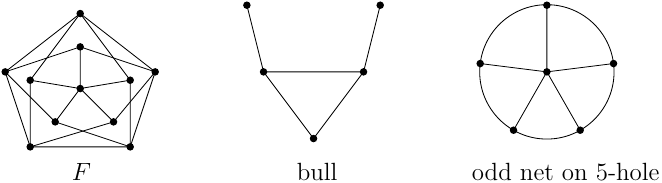}
		\end{center}
		\vskip -25pt
		\caption{Illustration of  some forbidden configurations.}
		\label{fig-1}
	\end{figure}

	Let $F$ denote the Gr\"{o}tzsch graph. A {\em bull} consists of a triangle with two disjoint pendant edges. Let $H$ be an odd hole or an odd antihole. An {\em odd net} is obtained from $H$ by adding a new vertex complete to $V(H)$; see Figure~\ref{fig-1}.
	
	A graph $G$ is {\em perfectly divisible} if, for every induced subgraph $H$ of $G$, $V(H)$ can be partitioned into $A$ and $B$ such that $H[A]$ is perfect and $\omega(H[B])<\omega(H)$. This notion was introduced by Ho\`{a}ng \cite{H22}. A simple induction on $\omega(G)$ shows that every perfectly divisible graph $G$ satisfies $\chi(G)\le\binom{\omega(G)+1}{2}$ \cite{CS19}.
	
	Let $G$ be a graph and let $w$ be a weight function on $V(G)$. Let $\omega_w(G)$ denote the maximum weight of a clique in $G$. A graph $G$ is {\em perfectly weight divisible} \cite{CS19,H22} if, for every positive integer weight function $w$ on $V(G)$ and every induced subgraph $H$ of $G$, there is a partition of $V(H)$ into $P$ and $W$ such that $H[P]$ is perfect and $\omega_w(H[W])<\omega_w(H)$.  If $G$ is not perfectly weight divisible but each of its proper induced subgraphs is, then $G$ is called {\em minimally non-perfectly weight divisible} (MNWD for short). Clearly, every perfectly weight divisible graph is perfectly divisible.
	
	Chudnovsky and Sivaraman \cite{CS19} proved that every (odd hole, bull)-free graph and every ($P_5$, bull)-free graph are perfectly divisible. This property was subsequently established for several other classes of bull-free graphs. In particular, Chen and Xu \cite{CX2025} proved that every ($P_7,C_5$, bull)-free graph is perfectly divisible. Clearly, $F$ is ($P_6$, bull)-free but not perfectly divisible, since $\chi(F)=4$ and $\omega(F)=2$. More generally, $F+K_n$ is ($P_6$, bull)-free, has clique number $n+2$, and is not perfectly divisible. Further results on perfect divisibility and bull-free graphs can be found in \cite{CW2025,DC2025,H26,KKS22}. Since every ($P_5$, bull)-free graph is perfectly divisible, it is natural to ask under what conditions a ($P_6$, bull)-free graph is perfectly divisible.
	
	We extend the above results in two directions. First, we extend the result of Chen and Xu \cite{CX2025} from ($P_7,C_5$, bull)-free graphs to ($P_8,C_5$, bull)-free graphs. Second, we characterize perfectly divisible ($P_6$, bull)-free graphs, thereby extending the result of Chudnovsky and Sivaraman \cite{CS19} on ($P_5$, bull)-free graphs.
	
	\begin{theorem}\label{P8}
		Let $G$ be a $(P_8,C_5, bull)$-free graph. Then $G$ is perfectly weight divisible, and hence $G$ is perfectly divisible.
	\end{theorem}
	
	\begin{theorem}\label{P6}
		Let $G$ be a $(P_6, bull)$-free graph. Then $G$ is perfectly divisible if and only if $G$ is $F$-free.
	\end{theorem}
	
	The remainder of the paper is organized as follows. In Section~\ref{notations}, we introduce some notation and preliminary results. In Sections~\ref{3} and~\ref{4}, we prove Theorems~\ref{P8} and~\ref{P6}, respectively, by analyzing the structure around a $7$-hole and a $5$-hole.
	
	\section{Notation and Preliminary Results}\label{notations}
	
	For $X\subseteq V(G)$, we use $G[X]$ to denote the subgraph of $G$ induced by $X$. Let $v\in V(G)$ and $X\subseteq V(G)$. We use $N_G(v)$ to denote the set of vertices adjacent to $v$, and let $M_G(v)=V(G)\setminus (N_G(v)\cup\{v\})$, $N_G(X)=\{u\in V(G)\setminus X\mid u$ has a neighbor in $X\}$, and $M_G(X)=V(G)\setminus (X\cup N_G(X))$. When there is no ambiguity, we omit the subscript $G$ and simply write $N(v)$, $M(v)$, $N(X)$ and $M(X)$.
	
	For a subset $A$ of $V(G)$ and a vertex $u\in V(G)\setminus A$, we say that $u$ is \emph{complete} to $A$ if $u$ is adjacent to every vertex of $A$, and that $u$ is \emph{anticomplete} to $A$ if $u$ is not adjacent to any vertex of $A$. For two disjoint subsets $A$ and $B$ of $V(G)$, $A$ is complete to $B$ if every vertex of $A$ is complete to $B$, and $A$ is anticomplete to $B$ if every vertex of $A$ is anticomplete to $B$.
	
	For $u,v\in V(G)$, we write $u\sim v$ if $uv\in E(G)$, and $u\not\sim v$ otherwise.
	
	Let $S\subseteq V(G)$ with $1<|S|<|V(G)|$. We say that $S$ is {\em a homogeneous set} of $G$ if every vertex in $V(G)\setminus S$ is either complete to $S$ or anticomplete to $S$. We use the following lemmas.
	
	\begin{lemma}\label{homogeneous} $($\em Theorem 3.6 of \cite{CS19}$)$
		No MNWD graph has a homogeneous set.
	\end{lemma}
	
	\begin{lemma}\label{MNPD}$($\em Lemma 6 of \cite{CWY}$)$
		Let $G$ be a bull-free MNWD graph. Then $G$ is a connected odd net-free graph, and for every vertex $v\in V(G)$, $G[M(v)]$ is imperfect and contains no odd antihole on at least seven vertices.
	\end{lemma}
	
	\begin{lemma}\label{P8C3}$($\em Theorem 2 of \cite{CS2018}$)$
		Let $G$ be a $(P_8,C_5,C_3)$-free graph. Then $\chi(G)\leq 3$.
	\end{lemma}
	
	\begin{lemma}\label{P6C3}$($\em Theorem 6 of \cite{RST2002}$)$
		Let $G$ be a $(P_6,C_3)$-free graph. Then $\chi(G)\leq 4$, and $G$ contains $F$ if $\chi(G)=4$.
	\end{lemma}
	
	\begin{lemma}\label{3-coloring}
		Let $G$ be a graph with $\chi(G)\leq 3$. Then $G$ is perfectly weight divisible.
	\end{lemma}
	\pf
	Let $w$ be a positive integer weight function on $V(G)$, and let $G'$ be an induced subgraph of $G$. Since $\chi(G')\leq 3$, there is a partition $(A,B)$ of $V(G')$ such that $G'[A]$ is bipartite and $B$ is an independent set. Let $X$ be the set of isolated vertices of $G'$ that belong to $B$. Then $G'[A\cup X]$ is perfect. Moreover, every vertex of $B\setminus X$ has a neighbor in $A$, and since all weights are positive, $\omega_w(G'[B\setminus X])<\omega_w(G')$. Thus $(A\cup X,B\setminus X)$ is the required partition. This completes the proof. \qed
	
	\section{Proof of Theorem~\ref{P8}}\label{3}
	
	In this section, we prove Theorem~\ref{P8}. We begin with the following lemma.
	
	\begin{lemma}\label{L-2}
		Let $G$ be a connected $(P_8,C_5, odd~net, bull)$-free graph with $\omega(G)\geq 3$, and let $v\in V(G)$ such that $v$ belongs to a maximum clique of $G$. If $G[M(v)]$ contains a $7$-hole, then $G$ has a homogeneous set.
	\end{lemma}
	
	\pf Let $C=v_1v_2\cdots v_5v_6v_7v_1$ be a 7-hole in $G[M(v)]$. Throughout this proof, subscripts are taken modulo 7. Let $i\in [7]$, define
	\begin{eqnarray*}
		X_i&=&\{u\in N(V(C))~|~N_{C}(u)=\{v_i\}\};\\
		Y_i&=&\{u\in N(V(C))~|~N_{C}(u)=\{v_i,v_{i+2}\}\};\\
		Z_i&=&\{u\in N(V(C))~|~N_{C}(u)=\{v_i,v_{i+1},v_{i+2}\}\};
	\end{eqnarray*}
	
	Let $X=\bigcup_{i=1}^7X_i$, $Y=\bigcup_{i=1}^7Y_i$ and $Z=\bigcup_{i=1}^7Z_i$. We next prove the following claim.
	\begin{claim}\label{N-2}
		$N(V(C))=X\cup Y\cup Z$.
	\end{claim}
	\pf
	It suffices to prove that $N(V(C))\subseteq X\cup Y\cup Z$. Let $u\in N(V(C))$. Without loss of generality, suppose that $u\sim v_1$. If $N(u)\cap V(C)=\{v_1\}$, then $u\in X_1$ and we are done. So we may assume that $u$ has a neighbor in $\{v_2,v_3,\cdots,v_7\}$.
	
	Suppose that $u$ is anticomplete to $\{v_2,v_7\}$. Then $u$ has a neighbor in $\{v_3,v_4,v_5,v_6\}$. Assume that $u$ is anticomplete to $\{v_3,v_6\}$. If $u\sim v_4$, then there is a 5-hole $uv_1v_2v_3v_4u$, a contradiction. So, $u\not\sim v_4$; similarly, $u\not\sim v_5$. Therefore, $u$ is not anticomplete to $\{v_3,v_6\}$ and by symmetry, we may assume that $u\sim v_3$.  If $u\sim v_6$, then $u$ has a neighbor in $\{v_4,v_5\}$ to avoid a 5-hole $uv_3v_4v_5v_6u$; and by symmetry, we may assume that $u\sim v_4$, but now $\{v_4,u,v_3,v_2,v_6\}$ induces a bull. So, $u\not\sim v_6$. To avoid a 5-hole $uv_5v_6v_7v_1u$, we have that $u\not\sim v_5$. Moreover, to avoid an induced bull on $\{u,v_3,v_4,v_2,v_5\}$, we have that $u\not\sim v_4$. Now, $N(u)\cap V(C)=\{v_1,v_3\}$, and then $u\in Y_1$. So, we may assume that $u$ is not anticomplete to $\{v_2,v_7\}$.
	
	Without loss of generality, $u\sim v_2$. To avoid an induced bull on $\{u,v_1,v_2,v_3,v_7\}$, we may assume that $u\sim v_3$. If $u\sim v_4$, then $u\sim v_5$ to avoid an induced bull on $\{v_3,u,v_4,v_1,v_5\}$; similarly, $u\sim v_7$. Since $G$ is odd net-free, $u$ is not complete to $V(C)$, and so $u\not\sim v_6$; but then $\{v_4,v_5,u,v_2,v_6\}$ induces a bull, a contradiction. So $u\not\sim v_4$, and similarly, $u\not\sim v_7$.
	
	If $u$ is complete to $\{v_5,v_6\}$, then $\{u,v_5,v_6,v_4,v_7\}$ induces a bull. So $u$ has at most one neighbor in $\{v_5,v_6\}$. If $u\sim v_5$, then there exists a 5-hole $uv_5v_6v_7v_1u$. So, $u\not\sim v_5$. Similarly, $u\not\sim v_6$. Then $N(u)\cap V(C)=\{v_1,v_2,v_3\}$; this implies $u\in Z_1$. This proves Claim~\ref{N-2}. \qed
	\begin{claim}\label{cc-4}
		If there exists $u\in N(v)$ such that $u\notin N(V(C))$, then $G$ has a homogeneous set.
	\end{claim}
	\pf
	Let $u\in N(v)$ such that $u\notin N(V(C))$. Then $u\in M(V(C))$. Since $G$ is connected, let $P=x_1x_2\cdots x_h$ be a shortest induced path connecting the edge $uv$ to $C$, where $h\geq 3$, $x_1\in\{u,v\}$, and $x_h\in V(C)$. Without loss of generality, $x_1=u$ and $x_h=v_1$. Since $x_{h-1}\in N(V(C))$, by symmetry, we may assume that $x_{h-1}\in X_1\cup Y_1\cup Z_1$ by Claim~\ref{N-2}. If $x_{h-1}\in X_1$, then there exists an induced $P_8=x_{h-2}x_{h-1}v_1v_2v_3v_4v_5v_6$.  If $x_{h-1}\in Z_1$, then $\{v_2,v_1,x_{h-1},v_7,x_{h-2}\}$ induces a bull. Therefore $x_{h-1}\in Y_1$.
	
	We have $h=3$, for otherwise there exists an induced $P_8=v_7v_6v_5v_4v_3x_{h-1}x_{h-2}x_{h-3}$. So  $x_2\sim v_1$ and $x_2\sim v_3$ as $x_{h-1}\in Y_1$. Consequently, to avoid an induced $P_8=v_7v_6v_5v_4v_3x_2x_1v$, we have that $v\sim x_2$. Let $S$ be the component of $G[N(x_2)\cap M(V(C))]$ which contains the edge $uv$. We claim that $V(S)$ is a homogeneous set of $G$.
	
	Suppose that there exists a vertex $t\in V(G)\setminus V(S)$ such that $t$ is neither complete nor anticomplete to $V(S)$. So there exist two vertices $t_1,t_2\in V(S)$ such that $tt_1t_2$ is an induced $P_3$. Clearly $t\notin V(C)$ as $t_1\in M(V(C))$.  We next show that
	\begin{equation}\label{e-1}
		t\not\sim x_2.
	\end{equation}
	
	Suppose, to the contrary, that $t\sim x_2$. By the choice of $S$, we have that $t\in N(V(C))$ as $t\notin V(S)$. Then $t$ has at least two neighbors in $V(C)$. Indeed, if $N(t)\cap V(C)=\{v_i\}$ for some $i\in[7]$, then $t_1tv_iv_{i+1}v_{i+2}v_{i+3}v_{i+4}v_{i+5}$ induces a $P_8$, a contradiction. Hence $t\in Y\cup Z$ by Claim~\ref{N-2}.
	
	Suppose that $t\not\sim v_1$. Moreover, $t$ is anticomplete to $V(C)\setminus\{v_2,v_3,v_7\}$ as otherwise $\{t,t_1,x_2,v_1,t'\}$ induces a bull, where $t'\in V(C)\setminus\{v_1,v_2,v_3,v_7\}$ is a neighbor of $t$, a contradiction. If $t\sim v_3$, then $\{t,x_2,v_3,t_2,v_2\}$ induces a bull whenever $t\not\sim v_2$, and $\{t,v_2,v_3,v_1,v_4\}$ induces a bull whenever $t\sim v_2$. Both cases yield contradictions. So $t\not\sim v_3$ and $t$ is anticomplete to $V(C)\setminus\{v_2,v_7\}$. By Claim~\ref{N-2}, $t\in Y_7$ as $t$ has at least two neighbors in $V(C)$. But then there exists an induced $P_8=v_6v_5v_4v_3v_2tt_1t_2$. So $t\sim v_1$, and by symmetry $t\sim v_3$.
	
	Now $t\in Y_1\cup Z_1$ by Claim~\ref{N-2}. But then $\{t,x_2,v_3,v_4,t_2
	\}$ induces a bull. This proves (\ref{e-1}).
	
	By (\ref{e-1}), we have that $t\sim v_3$ to avoid a bull on $\{t_2,t_1,x_2,t,v_3\}$. Similarly, $t\sim v_1$. By Claim~\ref{N-2}, $t\in Y_1\cup Z_1$. But then there exists an induced $P_8=v_7v_6v_5v_4v_3tt_1t_2$, a contradiction. Therefore $V(S)$ is a homogeneous set of $G$. This proves Claim~\ref{cc-4}. \qed
	
	By Claim~\ref{cc-4}, we may assume that
	\begin{equation}\label{eq-4}
		\mbox{	for every $u\in N(v)$, $u\in N(V(C))$.}
	\end{equation}
	
	\begin{claim}\label{c-3}
		For every $u\in N(v)$, $u\in Y$. And for any $x,y\in N(v)$ with $x\sim y$, $N_C(x)=N_C(y)$.
	\end{claim}
	\pf
	Let $u\in N(v)$. By (\ref{eq-4}), $u\in N(V(C))$. By Claim~\ref{N-2} and symmetry, we may assume that $u\in X_1\cup Y_1\cup Z_1$. If $u\in X_1$, then there is an induced $P_8=vuv_1v_2v_3v_4v_5v_6$. If $u\in Z_1$, then $\{v_2,v_3,u,v,v_4\}$ induces a bull. So $u\in Y_1\subseteq Y$.
	
	Let $x,y\in N(v)$ such that $x\sim y$. By symmetry, we may assume that $N_C(x)=\{v_1,v_3\}$ and $N_C(y)\in \{\{v_1,v_3\},\{v_2,v_4\}, \{v_3,v_5\}, \{v_4,v_6\}\}$. If $N_C(y)=\{v_2,v_4\}$, then $\{v,x,y,v_1,v_4\}$ induces a bull; if $N_C(y)=\{v_3,v_5\}$, then $\{v,x,y,v_1,v_5\}$ induces a bull; if $N_C(y)=\{v_4,v_6\}$, then $\{v,x,y,v_1,v_6\}$ induces a bull. So, $N_C(y)=N_C(x)=\{v_1,v_3\}$. This proves Claim~\ref{c-3}. \qed
	
	\medskip
	
	Recall that $v$ belongs to a maximum clique of $G$ and $\omega(G)\geq 3$. Let $S$ be a component of $N(v)$ which contains an edge. Clearly, $1<|V(S)|<|V(G)|$. We next prove that $V(S)$ is a homogeneous set of $G$. By Claim~\ref{c-3}, by symmetry, we may assume that $N(V(S))\cap V(C)=\{v_1,v_3\}$ as $S$ is connected.
	
	Clearly, $\{v,v_1,v_3\}$ is complete to $V(S)$, and $(N(v)\setminus V(S))\cup (V(C)\setminus\{v_1,v_3\})$ is anticomplete to $V(S)$. Suppose that $V(S)$ is not a homogeneous set of $G$. Then there exists a vertex $z$ which has a neighbor and a non-neighbor in $V(S)$. Clearly, $z\notin N(v)\cup\{v\}\cup V(C)$. Since $S$ is connected, there exists an edge $xy\in E(G[S])$ such that $z\sim x$ and $z\not\sim y$.
	
	We first claim that $z$ is anticomplete to $\{v_1,v_3\}$. Suppose that $z\sim v_1$. If $z\not\sim v_2$, then $\{x,z,v_1,v,v_2\}$ induces a bull, and hence $z\sim v_2$. Similarly, considering $\{x,z,v_1,v,v_7\}$, we obtain $z\sim v_7$. By Claim~\ref{N-2}, it follows that $N(z)\cap V(C)=\{v_7,v_1,v_2\}$. But then $\{z,v_7,v_1,v_6,y\}$ induces a bull, a contradiction. Thus $z\not\sim v_1$, and by symmetry $z\not\sim v_3$.
	
	Now, if $z\not\sim v_2$, then $\{x,y,v_1,z,v_2\}$ induces a bull. Hence $z\sim v_2$. Similarly, considering the triangles $xyv_3$ and $xyv_1$, we obtain $z\sim v_4$ and $z\sim v_7$, respectively. Therefore, $\{v_2,v_4,v_7\}\subseteq N(z)\cap V(C)$, while $z$ is anticomplete to $\{v_1,v_3\}$. Thus $z\in N(V(C))$, but $z\notin X\cup Y\cup Z$, contrary to Claim~\ref{N-2}. Therefore, $V(S)$ is a homogeneous set of $G$. This completes the proof of Lemma~\ref{L-2}. \qed
	
	\medskip
	
	We now prove Theorem~\ref{P8}.
	
	\noindent\textbf{{\em Proof of Theorem~\ref{P8}:}} On the contrary, let $G$ be a ($P_8,C_5$, bull)-free MNWD graph. By Lemma~\ref{MNPD}, we have that $G$ is a connected odd net-free graph, and for every vertex $v\in V(G)$, $G[M(v)]$ is imperfect and contains no odd antihole on at least seven vertices. We have $\omega(G)\geq 3$. Indeed, if $\omega(G)\leq 2$, then $G$ is triangle-free, and Lemma~\ref{P8C3} gives $\chi(G)\leq 3$. Lemma~\ref{3-coloring} would then imply that $G$ is perfectly weight divisible, a contradiction.
	
	Let $v\in V(G)$ such that $v$ belongs to a maximum clique of $G$. Since $G[M(v)]$ is imperfect and contains no odd antihole on at least seven vertices, Theorem~\ref{perfect} implies that $G[M(v)]$ contains an odd hole. As $G$ is $(P_8,C_5)$-free, this odd hole must be a 7-hole. Therefore, by Lemma~\ref{L-2}, $G$ has a homogeneous set, which contradicts Lemma~\ref{homogeneous}. This proves Theorem~\ref{P8}. \qed
	
	\section{Proof of Theorem~\ref{P6}}\label{4}
	
	In this section, we prove Theorem  \ref{P6}. We begin with the following lemma, whose proof follows the same general strategy as that of Lemma~\ref{L-2}.
	
	\begin{lemma}\label{L-1}
		Let $G$ be a connected $(P_6, odd~net, bull)$-free graph with $\omega(G)\geq 3$, and let $v\in V(G)$ such that $v$ belongs to a maximum clique of $G$. If $G[M(v)]$ contains a 5-hole, then $G$ has a homogeneous set.
	\end{lemma}
	\pf
	Let $C=v_1v_2\cdots v_5v_1$ be a 5-hole in $G[M(v)]$. Throughout this proof, subscripts are taken modulo 5. For
	$i\in [5]$, let
	\begin{eqnarray*}
		X_i&=&\{u\in N(V(C))~|~N_{C}(u)=\{v_i\}\};\\
		Y_i&=&\{u\in N(V(C))~|~N_{C}(u)=\{v_i,v_{i+2}\}\};\\
		Z_i&=&\{u\in N(V(C))~|~N_{C}(u)=\{v_i,v_{i+1},v_{i+2}\}\};\\
		W_i&=&\{u\in N(V(C))~|~N_{C}(u)=\{v_i,v_{i+1},v_{i+2},v_{i+3}\}\};
	\end{eqnarray*}
	
	Let $X=\bigcup_{i=1}^5X_i$, $Y=\bigcup_{i=1}^5Y_i$, $Z=\bigcup_{i=1}^5Z_i$ and $W=\bigcup_{i=1}^5W_i$. We next prove the following claim.
	\begin{claim}\label{N-1}
		$N(V(C))=X\cup Y\cup Z\cup W$.
	\end{claim}
	\pf
	It suffices to prove that $N(V(C))\subseteq X\cup Y\cup Z\cup W$. Let $u\in N(V(C))$. Without loss of generality, suppose that $u\sim v_1$. If $N(u)\cap V(C)=\{v_1\}$, then $u\in X_1$ and we are done. So, suppose that $u$ has a neighbor in $\{v_2,v_3,v_4,v_5\}$.
	
	If $u$ is anticomplete to $\{v_2,v_5\}$, then $u$ is not complete to $\{v_3,v_4\}$ as otherwise $\{u,v_3,v_4,v_2,v_5\}$ induces a bull. Hence, $u$ has exactly one neighbor in $\{v_3,v_4\}$. Now, $u\in Y_1$ if $u\sim v_3$ and $u\in Y_4$ if $u\sim v_4$. Therefore, we may assume that $u$ has a neighbor in $\{v_2,v_5\}$, and by symmetry assume that $u\sim v_2$.
	
	To avoid an induced bull on $\{u,v_1,v_2,v_3,v_5\}$, $u$ has a neighbor in $\{v_3,v_5\}$. By symmetry, we may assume that $u\sim v_3$. If $N(u)\cap V(C)=\{v_1,v_2,v_3\}$, then $u\in Z_1$ and we are done. So, we may assume that $u$ has a neighbor in $\{v_4,v_5\}$. Since $G$ is odd net-free, we have that $u$ is not complete to $V(C)$. Now, $u\in W_1$ if $u\sim v_4$ and $u\in W_5$ if $u\sim v_5$. This proves Claim~\ref{N-1}. \qed
	
	\medskip
	
	We next prove that
	\begin{claim}\label{cc-3}
		If there exists $u\in N(v)$ such that $u\notin N(V(C))$, then $G$ has a homogeneous set.
	\end{claim}
	\pf
	Let $u\in N(v)$ such that $u\notin N(V(C))$. Then $u\in M(V(C))$. Since $G$ is connected, let $P=x_1x_2\cdots x_h$ be a shortest induced path connecting the edge $uv$ to $C$, where $h\geq 3$, $x_1\in\{u,v\}$, and $x_h\in V(C)$. Without loss of generality, $x_1=u$ and $x_h=v_1$. Since $x_{h-1}\in N(V(C))$, by symmetry, we may assume that $x_{h-1}\in X_1\cup Y_1\cup Z_1\cup W_1$ by Claim~\ref{N-1}. If $x_{h-1}\in X_1$, then there exists an induced $P_6=x_{h-2}x_{h-1}v_1v_2v_3v_4$.  If $x_{h-1}\in Z_1\cup W_1$, then $\{v_2,v_1,x_{h-1},v_5,x_{h-2}\}$ induces a bull. Therefore $x_{h-1}\in Y_1$.
	
	We have $h=3$, for otherwise there exists an induced $P_6=v_5v_4v_3x_{h-1}x_{h-2}x_{h-3}$. So  $x_2\sim v_1$ and $x_2\sim v_3$ as $x_{h-1}\in Y_1$. Consequently, to avoid an induced $P_6=v_5v_4v_3x_2x_1v$, we have that $v\sim x_2$.
	Let $S$ be the component of $G[N(x_2)\cap M(V(C))]$ which contains the edge $uv$. We claim that $V(S)$ is a homogeneous set of $G$.
	
	Suppose, to the contrary, that there exists a vertex $t\in V(G)\setminus V(S)$ such that $t$ is neither complete nor anticomplete to $V(S)$. So there exist two vertices $t_1,t_2\in V(S)$ such that $tt_1t_2$ is an induced $P_3$. Clearly $t\notin V(C)$ as $t_1\in M(V(C))$. We next show that
	\begin{equation}\label{e-22}
		t\not\sim x_2.
	\end{equation}
	
	Suppose, to the contrary, that $t\sim x_2$. By the choice of $S$, we have that $t\in N(V(C))$ as $t\notin V(S)$. Clearly, $t$ has at least two neighbors in $V(C)$. Indeed, if $N(t)\cap V(C)=\{v_i\}$ for some $i\in[5]$, then $t_1tv_iv_{i+1}v_{i+2}v_{i+3}$ induces a $P_6$, a contradiction.
	
	Suppose that $t\not\sim v_1$. Then $t\not\sim v_4$ as otherwise $\{t,t_1,x_2,v_1,v_4\}$ induces a bull, a contradiction. Consequently, $t\not\sim v_3$ as otherwise $\{t,x_2,v_3,t_2,v_2\}$ induces a bull whenever $t\not\sim v_2$, and $\{t,v_2,v_3,v_1,v_4\}$ induces a bull whenever $t\sim v_2$. Both cases yield contradictions. So $t$ is anticomplete to $\{v_1,v_3,v_4\}$. By Claim~\ref{N-1}, $t\in Y_5$. But now there is an induced $P_6=v_4v_3v_2tt_1t_2$. Therefore $t\sim v_1$, and by symmetry $t\sim v_3$.
	
	Now, $t\sim v_2$ as otherwise $\{t,v_3,x_2,v_2,t_2\}$ induces a bull. By Claim~\ref{N-1}, by symmetry, we may assume that $t\in Z_1\cup W_1$. But then $\{v_2,v_1,t,v_5,t_1\}$ induces a bull. This proves (\ref{e-22}).
	
	By (\ref{e-22}), we have that $t\sim v_3$ to avoid a bull on $\{t_1,t_2,x_2,t,v_3\}$. Similarly, $t\sim v_1$. By Claim~\ref{N-1}, $t\in Y_1\cup W_3$. But then $v_5v_4v_3tt_1t_2$ is an induced $P_6$ if $t\in Y_1$, and $\{v_5,v_1,t,v_2,t_1\}$ induces a bull if $t\in W_3$. Both cases yield contradictions.  Therefore $V(S)$ is a homogeneous set of $G$. This proves Claim~\ref{cc-3}. \qed
	
	\medskip
	
	By Claim~\ref{cc-3}, we may assume that for every $u\in N(v)$,
	$$
	u\in N(V(C)).
	$$
	Consequently, we prove that
	\begin{equation}\label{e-2}
		\mbox{for every $u\in N(v)$, $u\in Y$.}
	\end{equation}
	
	By Claim~\ref{N-1}, without loss of generality, we assume that $u\in X_1\cup Y_1\cup Z_1\cup W_1$. If $u\in X_1$, then there is an induced $P_6=vuv_1v_2v_3v_4$. If $u\in Z_1$, then $\{v_2,u,v_3,v_4,v\}$ induces a bull. If $u\in W_1$, then $\{v_3,v_4,u,v,v_5\}$ induces a bull. So, $u\in Y_1$. This proves (\ref{e-2}).
	
	We next prove that
	\begin{equation}\label{e-3}
		\mbox{if $x,y\in N(v)$ such that $x\sim y$, then $N_C(x)=N_C(y)$.}
	\end{equation}
	
	Assume for contradiction that $N_C(x)\neq N_C(y)$. Without loss of generality, suppose that $N_C(x)=\{v_1,v_3\}$ by (\ref{e-2}). If $N_C(y)=\{v_2,v_4\}$, then $\{v,x,y,v_1,v_4\}$ induces a bull. If $N_C(y)=\{v_3,v_5\}$, then $\{v_3,x,y,v_2,v_5\}$ induces a bull. By symmetry, in all other cases, an induced bull also occurs. Hence, $N_C(x)=N_C(y)$. This proves (\ref{e-3}).
	
	Since $\omega(G)\geq 3$ and $v$ belongs to a maximum clique of $G$, there exist two adjacent vertices in $N(v)$. Let $S$ be a component of $N(v)$ which contains an edge. Clearly, $1<|V(S)|<|V(G)|$. We next prove that $V(S)$ is a homogeneous set of $G$. By (\ref{e-2}) and (\ref{e-3}), by symmetry, we may assume that $N(V(S))\cap V(C)=\{v_1,v_3\}$ as $S$ is connected.
	
	Clearly, $\{v,v_1,v_3\}$ is complete to $V(S)$, and $(N(v)\setminus V(S))\cup \{v_2,v_4,v_5\}$ is anticomplete to $V(S)$. Suppose that $V(S)$ is not a homogeneous set of $G$. There exists a vertex $z$ which has a neighbor and a non-neighbor in $V(S)$. Clearly, $z\notin N(v)\cup \{v\}\cup V(C)$. Since $S$ is connected, it follows that there exists an edge $xy\in E(G[S])$ such that $z\sim x$ and $z\not\sim y$.
	
	We first claim that $z$ is anticomplete to $\{v_1,v_3\}$. Suppose that $z\sim v_1$. If $z\not\sim v_2$, then $\{x,z,v_1,v,v_2\}$ induces a bull, and hence $z\sim v_2$. Similarly, considering $\{x,z,v_1,v,v_5\}$, we obtain $z\sim v_5$. By Claim~\ref{N-1}, $z\not\sim v_4$. But then $\{z,v_5,v_1,v_4,y\}$ induces a bull, a contradiction. Thus $z\not\sim v_1$, and by symmetry $z\not\sim v_3$.
	
	If $z\not\sim v_2$, then $\{x,y,v_1,z,v_2\}$ induces a bull, and hence $z\sim v_2$. Similarly, considering the triangles $xyv_3$ and $xyv_1$, we obtain $z\sim v_4$ and $z\sim v_5$, respectively. Therefore, $\{v_2,v_4,v_5\}\subseteq N(z)\cap V(C)$, while $z$ is anticomplete to $\{v_1,v_3\}$. Thus $z\in N(V(C))$, but $z\notin X\cup Y\cup Z\cup W$, contrary to Claim~\ref{N-1}. Therefore, $V(S)$ is a homogeneous set of $G$. This completes the proof of Lemma~\ref{L-1}. \qed
	
	\medskip
	
	We now prove Theorem~\ref{P6}.
	
	\noindent\textbf{{\em Proof of Theorem~\ref{P6}:}} Since perfect divisibility is hereditary and $F$ is not perfectly divisible, it suffices to prove that every ($P_6$, $F$, bull)-free graph is perfectly weight divisible.  Suppose not, and let $G$ be a ($P_6$, $F$, bull)-free MNWD graph. By Lemma~\ref{MNPD}, we have that $G$ is a connected odd net-free   graph, and for every vertex $v\in V(G)$, $G[M(v)]$ is imperfect and contains no odd antihole on at least seven vertices. We have $\omega(G)\geq 3$. Indeed, if $\omega(G)\leq 2$, then $G$ is triangle-free. By Lemma~\ref{P6C3}, $\chi(G)\leq 4$, and the equality case forces an induced copy of $F$. Since $G$ is $F$-free, $\chi(G)\leq 3$. Lemma~\ref{3-coloring} would then imply that $G$ is perfectly weight divisible, a contradiction.
	
	Let $v\in V(G)$ such that $v$ belongs to a maximum clique of $G$. Since $G[M(v)]$ is imperfect and contains no odd antihole on at least seven vertices, Theorem~\ref{perfect} implies that $G[M(v)]$ contains an odd hole. As $G$ is $P_6$-free, this odd hole must be a 5-hole. Therefore, by Lemma~\ref{L-1}, $G$ has a homogeneous set, which contradicts Lemma~\ref{homogeneous}. This completes the proof of Theorem~\ref{P6}. \qed
	
	{	\section*{Declarations}
		\begin{itemize}[nosep,leftmargin=*]
			\item \noindent{\bf Acknowledgements} The authors would like to thank Dr. Ran Chen for bringing this problem to their attention and for his valuable suggestions.
			\item \textbf{Conflict of interest}\quad The authors declare no conflict of interest.
			\item \textbf{Data availability statement}\quad This manuscript has no associated data.
	\end{itemize}}

\end{document}